\documentclass[12pt]{amsart}
\usepackage{amssymb,amsmath,amsthm,andy}

\newtheorem{thm}{Theorem}[section]

\newtheorem{cor}[thm]{Corollary}

\newtheorem{defn}[thm]{Definition}

\def\bar{\overline}
\def\wi{\widehat}
\def\n{\noindent}
\def\m{\medskip}
\renewcommand{\H}{\mathbb{H}}
\newcommand{\cLg}{\mathcal{L}_\gamma}
\newcommand{\trho}{\tilde\rho}

\begin{document} 

\begin{abstract} 
Let $\cLg = -1/4 \left( \sum_{j=1}^n(X_j^2+Y_j^2)+i\gamma T \right)$ where
$\gamma \in \C$, and 
$X_j, \ Y_j$ and $T$ are the left invariant vector fields
of the Heisenberg group structure for $\R^n \times \R^n \times \R$.
We explicitly compute the Fourier
transform (in the spatial variables) of the fundamental solution of the
Heat Equation $\p_s\rho = -\cLg\rho$. As a consequence, we have a simplified computation
of the Fourier transform of the fundamental solution of the $\Box_b$-heat equation on the
Heisenberg group and an explicit kernel of the heat equation associated to the weighted
$\dbar$-operator in $\C^n$ with weight $\exp(-\tau P(z_1,\dots,z_n))$ where 
$P(z_1,\dots,z_n) = \frac 12(|\Imm z_1|^2 + \cdots |\Imm z_n|^2)$ and $\tau\in\R$.
\end{abstract}

\title [Heat Equation on the Heisenberg Group] 
{A Simplified Calculation for the Fundamental Solution 
to the Heat Equation on the Heisenberg Group}

\author{Albert Boggess and Andrew Raich}

\address{
Department of Mathematics\\ Texas A\&M University\\ Mailstop 3368  \\ College Station, TX
77845-3368}

\subjclass[2000]{32W30, 33C45, 42C10}

\keywords{Heisenberg group, heat equation, fundamental solution, heat kernel, Kohn Laplacian}
\email{boggess@math.tamu.edu, araich@math.tamu.edu}
\maketitle



%


%


\setcounter{section}{-1}

\section{Introduction}\label{sec:intro}

The purpose of this note is to present a simplified calculation of the Fourier transform of
fundmental solution of the $\Box_b$-heat equation 
on the Heisenberg group. The Fourier transform of the fundamental solution 
has been computed
by a number of authors \cite{Ga77,Hu76,ChTi00,Tie06}. We use the approach of
\cite{ChTi00, Tie06} and compute the heat kernel using Hermite functions but
differ from the earlier approaches by working on a different, though biholomorphically
equivalent, version of the Heisenberg group. The simplification in the computation occurs
because the differential operators on this equivalent Heisenberg group take on a simpler form.
Moreover, in the proof of Theorem \ref{thm:fund soln transformed heat equation}, we reduce the $n$-dimensional heat
equation
to a 1-dimensional heat equation, and this technique would also be useful when analyzing the heat equation
on the nonisotropic Heisenberg group (e.g., see \cite{ChTi00}).
We actually use  the same version of the
Heisenberg group as Hulanicki \cite{Hu76}, but he computes the fundamental solution
of the heat equation associated to the sub-Laplacian and not the
Kohn Laplacian acting on $(0,q)$-forms.

A consequence of our fundamental solution computation is that we can explicitly compute
the heat kernel associated to the weighted $\dbar$-problem in $\C^n$ when the weight
is given by $\exp(-\tau P(z_1,\dots, z_n))$ where $\tau\in\R$ and
$P(z_1,\dots,z_n) = \frac 12(|\Imm z_1|^2 + \cdots |\Imm z_n|^2)$. When $n=1$ and
$p(z_1)$ is 
subharmonic, nonharmonic polynomial, 
the weighted $\dbar$-problem (with weight $\exp(-p(z_1))$) 
and explicit construction of Bergman and Szeg\"o
kernels and has been studied
by a number of authors in different contexts
(for example, see \cite{Christ91,Has94,Has95, Has98,FoSi91, Ber92}). In addition,
Raich has estimated the heat kernel and its derivatives
\cite{Rai06f, Rai06h,Rai07,Rai07p}.

%
%

\section{The Heisenberg Group and the $\Box_b$-heat equation}\label{sec:Heisenberg Group}

\begin{defn}\label{defn:Heisenberg group}

The Heisenberg group is the set  $\H^n=\R^n \times \R^n \times \R$ with the following group 
structure:
\[
g*g' = (x,y,t)* (x',y', t')= (x+x', y+y', t+t' +x\cdot y')
\]
where $(x,y,t),(x',y',t')\in \R^n\times\R^n\times\R$ 
and $\cdot$ denotes the standard dot product in $\R^n$.
\end{defn}
The left-invariant vector fields for this group structure are:
\[
X_j^g = \frac{\p}{\p x_j} + y_j \frac{\p}{\p t}
\ \ \textrm{and} \ \  Y_j^g = \frac{\p}{\p y_j}, \ \ 1 \leq j \leq n, \ \ \textrm{and}
\ \ T^g = \frac{\p}{\p t}.
\]
The Heisenberg group also can be identified with the following
hypersurface in $\C^{n+1}$: 
$H^n = \{(z_1,\dots,z_{n+1})\in\C^{n+1} : \Imm z_{n+1} = (1/2)\sum_{j=1}^n(\Imm z_j)^2\}$
where we identify $(z_1,\dots,z_n,t+i(1/2)\sum_{j=1}^n(\Imm z_j)^2)\in H^n$
with $(z_1,\dots,z_n,t)=(x_1,\dots,x_n,y_1,\dots,y_n,t)$ where $z_j = x_j+iy_j \in \C$.
With this identification, the left-invariant vector fields of type (0,1)
and (1,0), respectively are:
\[
\bar Z_j^g =(1/2)(X_j+iY_j)= \frac{\p}{\p\z_j} + \frac {y_j}2\frac{\p}{\p t}, 
\ \  \ \ Z_j^g = (1/2)(X_j-iY_j) =\frac{\p}{\p z_j} + \frac {y_j}2\frac{\p}{\p t}
\]
for $g = (x,y,t) \in \H^n$ and  $1\leq j\leq n$.



\m

\n \textbf{The Heat Equation.}  The Kohn Laplacian
$\Box_b$ acting on $(0,q)$-forms on $H^n \approx \H^n$
can be easily described in terms of these left-invariant
vector fields. Suppose $f = \sum_{J\in\I_q} f_J d\z_J$
is a $(0,q)$-form where $\I_q$ is the set of all 
increasing $q$-tuples $J=(j_1,\dots,j_q)$, $1\leq j_k\leq n$.
Then
\[
\Box_b f = \sum_{J\in\I_q} \cL_{n-2q} f_J\, d\z_J
\]
where 
\begin{equation}\label{eqn:L gamma}
\cLg = -\frac 14  \left( \sum_{j=1}^n (X_j^2 + Y_j^2) + i\gamma T \right) .
\end{equation}
See Stein (\cite{St93}, XIII \S2), for details on computing $\Box_b$.
For comparison, the box  operator ((or Laplacian) in Hulanicki (\cite{Hu76}) is 
$-\frac 12 \sum_{j=1}^n (X_j^2 + Y_j^2)$.

The Heat Equation is defined on $(0,q)$-forms $\rho$ 
on $\H^n$ with coefficient functions that 
depend on $s \in (0, \infty)$ and 
$(x,y,t) \in \H^n$. It is
\[
\frac{\p \rho}{\p s} = -\Box_b \rho
\]
(note that here, $s$ is the  ``time" variable and $t$ is 
a spatial variable). Since $\Box_b$ acts diagonally, 
we can restrict ourselves to a fixed component and look for a fundamental solution
$\rho$ that satisfies
\begin{equation}\label{eqn:heat 1}
\begin{cases} \displaystyle \frac{\p \rho}{\p s} = -\cL_\gamma \rho  \ \ \ \textrm{for} \ \ 
s >0, \ (x,y,t) \in \H^n \vspace{.1in} \\
\rho (s=0, x,y,t) = \delta_0 (x,y,t)
\end{cases}
\end{equation}
(i.e., the delta function at the origin in the spatial variables).

\m

\n \textbf{Fourier Transformed Variables.} We will use a Fourier 
transform in the spatial $(x,y,t)$ variables (i.e. {\em not}
the $s$-variable): let $(\alpha, \beta, \tau)$
be the transform variables corresponding to $(x,y, t)$, and define:
\[
\wi f (\alpha, \beta, \tau) 
= \int_{\H^n} f(x,y,t) \, e^{-i(\alpha\cdot x + \beta\cdot y + \tau t)} \, dx  \, dy \, d t.
\]
Our main result is the following:

\begin{thm}\label{thm:fund soln transformed heat equation}
For any $\gamma \in \C$, the spatial Fourier transform of the 
fundamental solution to the heat equation (\ref{eqn:heat 1}) is given by
\begin{equation}\label{eqn: fund soln}
\hat \rho^\gamma (s, \alpha, \beta, \tau)
= \frac{e^{-\gamma s \tau/4}} {(\cosh( s \tau/2))^{n/2} }
e^{-A ( |\alpha|^2+ |\beta|^2)/2 + iB \alpha\cdot \beta}.
\end{equation}
where
\[
A= \frac{\sinh(s \tau/2)}{\tau \cosh (s \tau/2)} , \ \ \ \ 
B= \frac{2 \sinh^2 (s \tau/4)}{\tau \cosh (s \tau/2)}.
\]
\end{thm}
Note that $\gamma$ may be any complex number, but $\gamma = n-2q$ is the value
where $\cLg$ corresponds to
$\Box_b$ on $(0,q)$-forms.

We also seek the fundamental solution to the heat equation
associated to the weighted $\dbar$ operator in $(s,x,y)$-space.
Given a function $f$ on $\R^n\times\R^n\times\R$, let 
\[
\tilde f_\tau(x,y) = \int_\R e^{-i\tau t}f(x,y,t)\, dt
\]
be the partial Fourier transform in $t$.
Define
\[
\bar L_j = \frac{\p}{\p\z_j} + \frac i2 y_j\tau
= \frac 12 (\frac{\p}{\p x_j} +i\frac{\p}{\p y_j} +iy\tau ), \ \ \ 
L_j = \frac{\p}{\p z_j} + \frac i2 y_j\tau
= \frac 12 (\frac{\p}{\p x_j} - i\frac{\p}{\p y_j} +iy\tau ).
\]
Note that these operators are just the Fourier transform of 
$\ZZ_j$ and $Z_j$ in the $t$-direction. If $\triangle_{x,y}$ is the Laplacian
in both the $x$ and $y$ variables, 
the partial $t$-Fourier transform of $\cLg$ is
\[
\tilde \cLg = -\frac 14 \big( \triangle_{x,y} + 2i\tau\, y\cdot\nabla_x
- (\tau^2 y\cdot y + \gamma\tau)\big).
\]
The operator $\tilde\cLg$ acts on functions, but it can 
be extended to $(0,q)$-forms by acting on each component function
of the form. If $\gamma=n-2q$, then $\tilde\cLg$ is the higher dimensional analog of the
$\Boxtp$-operator from \cite{Rai06h,Rai07,Rai07p} 
associated to the weighted $\dbar$ operator in $\C^n$
with weight $\exp(-\tau P(z_1,\dots,z_n))$ where 
$P(z_1,\dots,z_n) = \frac 12(|\Imm z_1|^2 + \cdots |\Imm z_n|^2)$ and $\tau\in\R$.  As a corollary to our main theorem,
we compute the fundamental solution to the heat operator associated
to this weighted $\dbar$.

\begin{cor}\label{cor:trho calculation}
For any $\gamma \in \C$, $\tau \in \R$, the function 
\[
\trho_\tau^\gamma (s, x,y) 
= \frac{ e^{-\gamma s \tau/4}}{(2\pi)^{n} (\cosh( s \tau/2))^{n/2} (A^2+B^2)^{n/2} } 
e^{{-\frac A{2(A^2+B^2)} (|x|^2+|y|^2)} -{i\frac B{A^2+B^2}  x\cdot y} }. 
\]
is the fundamental solution to the weighted $\dbar$ heat equation:
$(\frac \partial {\partial s} + \tilde \cLg )\trho_\tau^\gamma (s, x,y ) =0$
with $\trho_\tau^\gamma (s=0, x,y) = \delta_{(0,0)} (x,y)$.
\end{cor}

Finally, we use $\trho_\tau^\gamma$ to derive the heat kernel,
as studied in \cite{Rai06h,Rai07,Rai07p,NaSt01h}.

\begin{cor}\label{cor: heat kernel}
For any $\gamma \in \C$, $\tau \in \R$, let
\[
H_\tau^\gamma(s,x', y', x,y)=
\frac{\tau^n e^{-\gamma s \tau/4}} {(4\pi)^n \sinh^n(s \tau/4)}
e^{ -\frac\tau 4 \coth(s \tau/4) 
\left( |x-x'|^2+|y-y'|^2 \right) -{i \frac\tau2}
(x-x')\cdot(y+y')}.
\]
Then $H_\tau^\gamma$ is the heat kernel which satisfies the following property:
if $f \in L^2 (\C)$, then
\[
H_\tau^\gamma [f](s, x,y) = 
\int_{\R^n\times\R^n} H_\tau^\gamma(s, x, y, x', y') f(x',y') \, dx' dy' 
\]
is a solution to the following initial value problem for the heat equation: 
\begin{equation}\label{eqn:heat kernel equation}
\begin{cases} \Big( \displaystyle \frac{\p}{\p s}+\tilde\cLg\Big)  H_\tau^\gamma [f] =  0 
   \vspace*{.1in}\\
 H_\tau^\gamma [f]  (s=0, x,y)= f(x,y).
\end{cases}
\end{equation}

\end{cor}
Note that $H_\tau^\gamma$ is conjugate symmetric in $z=x+iy$ and $z'=x'+iy'$
(i.e. switching $z$ with $z'$ results in a conjugate).

%
%
\section{Proof of Theorem \ref{thm:fund soln transformed heat equation}}

It is easy to verify the following calculations. Recall that
$\ \wi{}\ $ refers to spatial Fourier transform.
\begin{eqnarray*}
\wi{X_j^2 f} (\alpha, \beta, \tau) &=& (- \alpha_j^2 -2i \alpha_j \tau \frac{\p}{\p \beta_j} 
+ \tau^2 \frac{\p^2}{\p \beta_j^2} ) \wi f \\
\wi{Y_j^2 f } (\alpha, \beta, \tau) &=& - \beta_j^2 \wi f \\
\wi{T f} (\alpha, \beta, \tau) &=& i \tau \wi f.
\end{eqnarray*}
We first reduce the problem down to dimension one.
Define $\hat \rho^{\gamma, 1}$ by the same
formula as given in (\ref{eqn: fund soln}), but for dimension one
(i.e. $n=1$ and $\alpha, \ \beta \in \R$).
From (\ref{eqn: fund soln}), note that 
\begin{equation}
\label{product}
\hat \rho^\gamma (s, \alpha, \beta, \tau) = \prod_{j=1}^n \hat \rho^{\gamma/n, 1} (s, \alpha_j, 
\beta_j,
\tau), \ \ \alpha =(\alpha_1, \dots, \alpha_n), \ \beta = (\beta_1, \dots , \beta_n) \in \R^n
\end{equation}
(note the $\gamma$ on the left and the $\gamma/n$ on the right).
Once we show that $\rho^{\gamma, 1}$ satisfies the 
transformed heat-equation in dimension one, i.e., 
\begin{equation}\label{eqn:transformed 1-d equation}
\left( \frac \partial {\partial s} - (1/4)(\wi{X^2}+ \wi{Y^2} +i \gamma \wi{T}) \right)
\{ \wi{\rho}^{\gamma, 1} (s, \cdot , \cdot) \} =0
\end{equation}
with initial condition $\hat \rho^{\gamma, 1}  = 1$
(the Fourier transform of the delta function), 
then by using (\ref{product}),  it is an easy exercise to show that $\wi{\rho}^\gamma$ 
in dimension $n$ satisfies
Theorem \ref{thm:fund soln transformed heat equation}.
  
From now on, we assume the dimension $n$ is one and 
so $x,y, \alpha$ and $\beta$ are all real variables. Also, $\gamma$ will be suppressed
as a superscript.
Define
\begin{equation}
\label{rho}
u(s, \alpha , \beta, \tau) = \wi{\rho} (s, \alpha, \beta, \tau ) 
e^{- i {\frac{\alpha \beta}\tau}}.
\end{equation}
Then, the following equations are easily verified
\begin{eqnarray}
\label{eqn:u, s=0 condition}
u(s=0, \alpha, \beta, \tau) &=& e^{- i {\frac{\alpha \beta}\tau}} \\ 
\label{heattrans}
\frac{\p u}{\p s} &=& (1/4) ( \tau^2 \frac {\partial^2}{\partial \beta^2} - \beta^2 - 
\gamma\tau) u .
\end{eqnarray}
The first equation follows from the fact that the Fourier transform of the delta
function is the constant one. The second equation follows from 
the heat equation for $\wi\rho$ (from (\ref{eqn:transformed 1-d equation})) 
and the above formulas for the transformed
differential operators $\wi X, \wi Y$ and $\wi T$. We will refer to the above differential
equation as the {\em transformed Heat equation}.

\m

\n \textbf{Solution of Heat Equation Using Hermite Special Functions.}
For $m=0, 1,2 \dots$ and $x\in\R$, let
\begin{eqnarray*}
\psi_m (x) &=& \frac{(-1)^m} 
{\sqrt{2^m m! \sqrt{\pi}}} e^{x^2/2} \frac{d^m}{dx^m}\{e^{-x^2}\}. \\
\end{eqnarray*}
For $\tau \in \R$, let
\begin{align*}
\Psi^\tau_m(x) &= |\tau|^{-1/4} \psi_{m}(x/\sqrt{ |\tau|}).
\end{align*}
It is a fact that $\psi_m$ and hence $\Psi^\tau_m$
form an orthonormal system for $L^2 (\R)$ (see \cite{Th93}, pg.1-7). 
It is also a fact (again see \cite{Th93}, (1.1.28)) that
\[
\psi_m''(x) = x^2 \psi_m (x) - (2m+1) \psi_m (x).
\]

We first assume that $\tau >0$ and later indicate the minor changes
needed  in the case that $\tau \leq 0$.
Replacing $x$ by $\beta/\sqrt{\tau}$ in the previous equation yields:
\begin{equation}
\label{eigenfunc}
(\tau^2 \frac {\partial^2}{\partial \beta^2}  - \beta^2 -\gamma\tau ) \{ \Psi_m^\tau \}
= - (2 m + 1 + \gamma)\tau \Psi^\tau_{m} (\beta). 
\end{equation}
In other words, $\Psi_m^\tau$ is an eigenfunction of the differential operator
on the right side of (\ref{heattrans}) with eigenvalue $-(2m +1+\gamma)\tau$.

Since $\{\Psi_m^\tau\}$ are an orthonormal basis for $L^2(\R)$, 
$u$ can be expressed as 
\[
u(s, \alpha, \beta, \tau) = \sum_{m = 0}^\infty
a_{m}( \alpha, \tau) e^{-\frac 14(2m +1+\gamma)s\tau} \Psi^\tau_{m} (\beta)
\]
where $a_m (\alpha, \tau)$ will be determined later.
Differentiating this with respect to $s$ and using (\ref{eigenfunc}) gives
\begin{align*}
\frac{\p}{\p s} u (s, \alpha, \beta, \tau) 
&= \sum_{m = 0}^\infty a_{m} ( \alpha, \tau) 
e^{-\frac 14(2m +1+\gamma)s\tau} (-\frac 14(2m +1+\gamma))\tau
 \Psi^\tau_{m} ( \beta) \\
&= \frac 14 \left( \tau^2 \frac {\partial^2}{\partial \beta^2} - \beta^2 - \gamma\tau \right) 
\{ u(t, \alpha, \beta, \tau) \}. 
\end{align*}
So, $u$ satisfies the transformed Heat equation (\ref{heattrans}). 
To satisfy the initial condition (\ref{eqn:u, s=0 condition}), we must have
\[
e^{-i \alpha \beta/\tau} = u(s=0, \alpha, \beta, \tau) = \sum_{m = 0}^\infty
a_{m}( \alpha, \tau) \Psi^\tau_{m} (\beta).
\]
Using the fact that the $\Psi^\tau_{m} (\beta) $ is an orthonormal system, we have
\begin{align*}
a_{m} (\alpha, \tau) 
&= \int_{\R} e^{-i \alpha \beta/\tau} \Psi^\tau_{m} ( \beta) \, d \beta 
= \tau^{1/4} \int_{\R} e^{-i \frac{\alpha}{\sqrt{\tau}}\beta} 
\psi_{m}(\beta)\, d\beta.
\end{align*}
The integral on the right is just the Fourier transform of $\psi_m$ 
at the point $\alpha/\sqrt{\tau}$.
From Thangavelu (\cite{Th93}, Lemma 1.1.3),  
the Fourier transform of $\psi_m$ equals $\psi_m$ up to a
constant factor of $(-i)^m \sqrt{2 \pi}$.
Therefore, 
\[
a_{m} (\alpha, \tau) =(-i)^{m} (2 \pi)^{1/2} \tau^{1/4} \psi_{m}( \alpha/\sqrt{\tau}). 
\]
Substituting this value of $a_{m}$ into the expression for $u$
and rearranging gives:
\[
u(s, \alpha, \beta, \tau)= (2 \pi)^{1/2} e^{-\frac 14(1+\gamma)s \tau}
\sum_{m = 0}^\infty
(-i)^{m} \psi_{m} (\frac{\alpha}{\sqrt{\tau}}) \psi_{m} (\frac{\beta}{\sqrt{\tau}})
e^{-\frac 12 m s\tau}.
\]
Now solving for $\hat \rho$ (see equation (\ref{rho})) yields
\[
\hat \rho(s, \alpha, \beta, \tau) = e^{i \alpha \beta/ \tau} u (s, \alpha, \beta, \tau) 
= (2 \pi)^{1/2} e^{-\frac 14(1+\gamma)s \tau}
\sum_{m = 0}^\infty
(-i)^{m} \psi_{m} (\tfrac{\alpha}{\sqrt{\tau}}) \psi_{m} (\tfrac{\beta}{\sqrt{\tau}})
e^{-\frac 12 m s\tau}e^{ i \alpha \beta/\tau}.
\]
Now let $S=e^{-s \tau/2}, \ x= \alpha/\sqrt{\tau}, \ y= \beta/\sqrt{\tau}$.
Since $|iS|<1$, we obtain (see \cite{Th93}, (1.1.36))
\begin{align*}
\hat \rho(s, \alpha, \beta, \tau) &= (2 \pi)^{1/2} S^{\frac 12(1+\gamma)} 
\left( \sum_{m = 0}^\infty (-iS)^{m} 
\psi_m (x) \psi_m (y)  \right) e^{ix y} \\
&= \frac{\sqrt2 S^{\frac 12(1+\gamma)}} {(1+S^2)^{1/2}}
e^{-\frac 12 \frac{1-S^2}{1+S^2}(x^2+y^2)} e^{ix y (\frac{-2S}{1+S^2} + 1)}. 
\end{align*}
Now substituting in for $S$, $x$ and $y$, a short calculation finishes the proof for $\tau>0$.
Note that $\hat \rho (s=0, \alpha, \beta, \tau) = 1$
(the Fourier transform of the delta function at the origin).

When $\tau=0$, the solution in (\ref{eqn: fund soln}) becomes
$\wi{\rho} (s, \alpha, \beta) = e^{-s(\alpha^2+\beta^2)/4}$
which is easily shown to satisfy (\ref{eqn:transformed 1-d equation}).

If $\tau <0$, then $\tau$ is replaced by $|\tau|$ on the right side of 
(\ref{eigenfunc}), which  slightly changes the subsequent calculations.
However the formula for the solution
given Theorem \ref{thm:fund soln transformed heat equation} remains valid
for $\tau <0$.

%
%

\section{Proof of the Corollaries}
\begin{proof}[Proof. (Corollary \ref{cor:trho calculation})]
Again, we assume the dimension is $n=1$. The fundamental
solution to this heat operator must satisfy
\[
\frac{\p}{\p s} \trho_\tau (s, x,y) +\tilde \cLg \trho_\tau=0 
\]
with the initial condition $\trho_\tau (s=0, x,y) = \delta_0(x,y)$.
Now since $\hat \rho$ is the Fourier transform of the fundamental
solution to the original Heat operator, clearly $\trho_\tau$
can be obtained by taking the inverse Fourier transform
of $\hat \rho$ in the $\alpha$, $\beta$ variables. This is a standard
calculation involving Gaussian integrals and will be left to the reader.
\end{proof}

\begin{proof}[Proof. (Corollary \ref{cor: heat kernel})]
If $L_j$ and $\bar L_j$, $1\leq j\leq n$, had constant coefficients then the 
heat kernel would just be $\trho_\tau (s, x-x',y-y')$ -- 
an ordinary convolution.
However, we must multiply by a ``twist'' factor $e^{- i \tau (x-x')\cdot y'}$
to account for the fact that $L_j$ and $\overline{L_j}$ have variable
coefficients. 
Let
\begin{equation}
\label{eqn:fund soln}
H_\tau(s, x, y, x',y', \tau) = \trho_\tau (s, x-x',y-y') e^{- i \tau (x-x')\cdot y'}.
\end{equation}
Note that $H_\tau (f)$ satisfies the initial condition 
given in (\ref{eqn:heat kernel equation}) in view of the
initial condition satisfied by $\tilde \rho_\tau$ and 
noting that the twist term 
is $1$ at $x'=x$. Showing that $H_\tau$ satisfies
the heat equation in the $s, \  x, \ y$ variables
is a short calculation that uses the equation
\[
\left( \frac{\p}{\p s} - \frac14 \Big(\triangle_{x,y} 
+2i \tau (y-y') \cdot\nabla_x
-(\tau^2 (y-y')\cdot(y-y') + \gamma \tau) \Big) \right) \{ \trho_\tau (s, x-x', y-y') \}
=0.
\]
which is just the equation $(\frac{\p}{\p s} +\tilde\cLg) \trho_\tau =0$
at the point $(s, x-x', y-y')$.
\m

\n \textbf{Simplification of the Formula for $H_\tau$}. Note that the 
coefficient of the imaginary part of the exponent of $\trho_\tau$ is
\[
\frac{- B}{A^2+B^2} \ \ \ 
\textrm{where} \ \ \ A= \frac{\sinh(s \tau/2)} {\tau \cosh (s \tau/2)} , \ \ 
B= \frac{2 \sinh^2 (s \tau/4)} {\tau \cosh (s \tau/2)}.
\]
An easy calculation with $\cosh$ and $\sinh$ identities shows that
\[
\frac{B}{A^2+B^2} = \frac \tau 2 \ \ \ \textrm{and}\ \ \ 
\frac{A}{B} = \frac{\cosh(s\tau/4)}{\sinh(s\tau/4)}.
\]
Consequently, the fundamental solution $H_\tau$, from 
\eqref{eqn:fund soln} and Corollary \ref{cor:trho calculation},
can be rewritten
\[
H_\tau(s,x', y', x,y)=
\frac{\tau^n e^{-\gamma s \tau/4}} {(4\pi)^n \sinh^n(s \tau/4)}
e^{ -\frac\tau 4 \coth(s \tau/4) 
\left( |x-x'|^2+|y-y'|^2 \right) -{i \frac\tau2}
(x-x')\cdot(y+y')}.
\]
\end{proof}
%
\bibliographystyle{alpha}
\bibliography{mybib}

\end{document}